\documentclass[12pt]{article}

\usepackage{fullpage}

\usepackage{amsfonts}

\usepackage{graphicx}

\usepackage{mathrsfs}

\usepackage{amsmath}

\usepackage{amssymb}

\usepackage{float}
\usepackage{color}
\usepackage{booktabs}
\usepackage[authoryear]{natbib}

\usepackage{geometry}
\newtheorem{theorem}{Theorem}

\numberwithin{equation}{section}

\def\cp{\mathop{\rightarrow}\limits^{p}}
\def\cd{\mathop{\rightarrow}\limits^{d}}

\allowdisplaybreaks[4]
\def\var{\mathrm {Var}}

\def\cov{\mathrm {Cov}}

\def\O{\bm\Omega}
\def\S{{\bf S}}

\def\bmv{\bm \varepsilon}

\def\A{{\mathcal A}}

\def\z{{\bf z}}

\def\I{{\bf{I}}}

\newcommand{\bm}{\boldsymbol}

\def\var{\mathrm {Var}}

\def\cov{\mathrm {Cov}}

\def\U{\bm U}

\def\u{\bm u}
\def\tr{{\rm tr}}
\def\bms{{\bf \Sigma}}

\def\beps{\boldsymbol e}

\def\z{{\bf z}}

\def\bSig{{\bf \Sigma}}

\def\A{{\bf A}}

\def\I{{\bf I}}

\def\0{{\bf 0}}
\def\1{{\bf 1}}

\title
{\bf Spatial-sign based High Dimensional White Noises Test}
\author{Ping Zhao, Dachuan Chen and Zhaojun Wang\\
Nankai University}
\date{}

\begin{document}

\maketitle

\begin{abstract}
A spatial-sign based test procedure is proposed for high dimensional white noise test in this paper. We establish the limit null distribution and give the asymptotical relative efficient of our test with respect to the test proposed by \citet{flm2022} under some special alternative hypothesis. Simulation studies also demonstrate the efficiency and robustness of our test for heavy-tailed distributions.
\end{abstract}

\noindent%
{\it Keywords:}  High-dimensional data, Spatial-Sign, White noise test

\section{Introduction}
In this paper, we consider testing for white noise or serial correlation, which is a fundamental problem in statistical inference. For univariate time series, the famous Box-Pierce portmanteau test and its variations are very popular due to their convenience in practical application \citep{Li2004,Lutkepohl2005}. Many efforts have been devoted to extending those methods for testing multivariate time series, such as \citet{Hosking1980,1981Distribution}. Recently, high dimensional time series data frequently appear in many applications, including finance and econometrics, biological and environmental research, etc, where the dimension of the time series are comparable or even larger than the observed length of the time series. In this case, the above traditional white noise tests can not directly apply for high dimensional data.

Recently, there are two types of omnibus tests proposed to deal with high dimensional white noise test. One is the max-type test. \citet{asw066} proposed a test statistic by using the maximum absolute auto-correlation and cross-correlations of the component series. \citet{2020Testing} proposed a rank-based max-type test by using the Spearman's rank correlation. \citet{cf2022} extend \citet{2020Testing}'s work to other types rank-based correlations, such as Kendall's tau correlation and Hoeffding’s D statistic, etc. As known to all, the max-type tests perform well for the sparse alternatives where only a few auto-correlations are nonzero and large, but perform less powerful for the dense alternatives where there are many small nonzero auto-correlations. So researchers proposed sum-type tests for the high dimensional white noise test. \citet{2019On} proposed a test statistic by using the sum of the squared singular values of several lagged sample autocovariance matrices. \citet{flm2022} proposed a new sum-type test statistics by excluding some terms in the test statistic proposed by \citet{2019On} and show that it has a better size performance. However, the above two sum-type tests are all based on the independent component model, which only allows the underling distribution of the time series is light tailed. Unfortunately, the assumption of light tailed distribution may be no appropriate for many applications, such as stock security returns. Thus, we need to construct a robust high dimensional white noise test procedure for the heavy tailed distributions.

The classic spatial sign based procedures are very robust and efficient in traditional multivariate analysis, see \citet{Oja2010Multivariate} for an overview. Recently, many literatures show that the spatial sign based procedures also perform very well in high dimensional settings. \citet{wpl}, \citet{fs2016}, \citet{f2021} proposed some spatial-sign based test procedures for the high dimensional one sample location problem. \citet{fzw2016}, \citet{h2022} consider high dimensional two sample location problem. \citet{zpfw2014}, \citet{fl2017} also extend the spatial-sign based method to the high dimensional sphericity test. Some spatial-sign based test procedures for high dimensional alpha test in factor pricing model are proposed by \citet{liu2023}, \citet{z2022}, \citet{z2023}. In an important work, \citet{paindaveine2016high} proposed a spatial-sign based test for i.i.d-ness against serial dependence. However, they assume the random vectors have independent spherical directions, which is too limited in applications. In practice, there are always some correlation between the random vectors. So we propose a new spatial sign based test procedure for the high dimensional white noise test in this article. Under the elliptical symmetric distribution assumption, we establish the asymptotical normality of the proposed test statistic under the null hypothesis and a special alternative hypothesis. We also show that the asymptotical relative efficiency of our method with respect to the test proposed by \citet{flm2022} is equivalent to the corresponding asymptotical relative efficiency of spatial-sign based method with respect to the least-square based procedures in high dimensional settings \citep{wpl,fs2016,liu2023,z2023}. Simulation studies also demonstrate the superiority of our method for heavy-tailed distributions.

This paper is organized as follows. In section 2,  we introduce our proposed spatial-sign based test procedure for high dimensional white noise test and establish the theoretical results. Simulation studies are showed in Section 3. All the technical details are collected in Section 4.

\section{Test Procedure}
Let $\bmv_1,\cdots,\bmv_n$ be a $p$-dimensional weakly stationary time series with mean zero.
We consider the following testing problem:
\begin{align}
H_0:\left\{\varepsilon_t\right\} is~ white~ noise~ v.s. ~H_1:\left\{\varepsilon_t\right\} ~is~ not~ white~ noise,
\end{align}
where the dimension of time series $p$ is comparable to or even greater than the sample size
$n$. Under the null hypothesis, $E(\bmv_t\bmv_{t+k}^\top)=\bm 0$. So \citet{2019On} proposed the following test statistic
\begin{align*}
G_H=\sum_{h=1}^H\tr(\hat{\S}_{h}^\top \hat{\S}_h), ~\hat{\S}_h=\frac{1}{n}\sum_{t=h+1}^n\bmv_t\bmv_{t-h}^\top.
\end{align*}
They established the asymptotical normality of $G_H$ by random matrix theory. \citet{flm2022} remove the diagonal elements $\bmv_i^\top \bmv_i\bmv_{i+h}^\top \bmv_{i+h}$ from the summation and proposed the following test statistic
\begin{align*}
T_{FLM}=\frac{1}{n(n-1)}\sum_{h=1}^H \underset{s\not=t}{\sum\sum}\bmv_t^\top \bmv_s\bmv_{t+h}^\top \bmv_{s+h}.
\end{align*}
They also established the asymptotical normality of $T_{FLM}$ by martingale central limit theorem. Both the above two tests need the independent component model assumption, i.e. $\bmv_t=\S^{1/2}\z_t$ and $\z_t=(z_{t1},\cdots,z_{tp})^\top$ is a sequence of $p$-dimensional independent random vectors with independent components. However,  A common drawback of the independent component model is their inability to handle many well-known heavy-tailed distributions, such as the multivariate Student $t$ and the mixture of multivariate normal distributions. Thus, we need to propose a robust and efficient test procedure for heavy-tailed distributions.

Under the assumption that $\bmv_t$ have independent spherical directions, \citet{paindaveine2016high} proposed a standardization test statistic
\begin{align}\label{tpv}
T_{PV}=\frac{\sqrt{2p^2}}{\sqrt{H}}\sum_{h=1}^H\frac{1}{n-h}\sum_{h+1\le s<t\le n}\U_{s-h}^\top \U_{t-h}\U_s^\top \U_t
\end{align}
where $\U_t=U(\bmv_t)$ and $U(\bm x)=\frac{\bm x}{||\bm x||}I(\bm x\not=0)$. They show that $T_{PV}\cd N(0,1)$ as $n,p\to \infty$ under the null hypothesis. However, the assumption of independent spherical directions always do not hold in practice. In addition, they do not give the power function of $T_{PV}$ under the alternative hypothesis. So we need to establish the theoretical results of the spatial-sign based test statistic under more general scatter matrix assumption.

We consider the following test statistics
\begin{align}
T_S=\sum_{h=1}^H\frac{1}{n-h}\sum_{h+1\le s<t\le n}\U_{s-h}^\top \U_{t-h}\U_s^\top \U_t
\end{align}
which is mimic to the test statistic (\ref{tpv}). Next, we will show that the asymptotic variance of $T_S$ under the null hypothesis is $\frac{H}{2}\tr^2(\O^2)$ and $\O=E(\U_t\U_t^\top)$, which is equal to $\frac{H}{2}p^{2}$ if $\bmv_t$ have independent spherical directions.

We need the following conditions:
\begin{itemize}
\item[(C1)] {\it (Error Distribution)} The error vectors ${\bmv}_1,\ldots,{\bmv}_n$ are
i.i.d. from the $p$-variate mean zero elliptical distribution with probability density function:
$$ \label{elliptical}
\mbox{det}({\bf \Xi})^{-1/2}g(\|
 {\bf \Xi}^{-1/2}{\bmv}\|),~{\bmv}\in \mathbb{R}^p.
$$
where ${\bf \Xi}$ is a positive definite scatter matrix.
\item[(C2)] {\it (Covariance Matrix)} $\tr(\bSig^4)=o(\tr^2(\bSig^2))$ and $\frac{\tr^4(\bSig)}{\tr^2\left(\bSig^2\right)} \exp \left\{-\frac{\tr^2(\bSig)}{128 N \lambda_{\max }^2(\bSig)}\right\}\to 0$ where  $\bSig=\cov(\beps_t)\doteq (\sigma_{ij})_{1\le i,j\le N}$ and $\lambda_{\max}(\bSig)$ is the largest eigenvalue of $\bSig$.
\end{itemize}
Under the condition (C1), $\bmv_i$ can be decomposed as ${\bf \Xi}^{1/2}R_i u_i$ where $u_i$ is a random vector uniformly distributed on the unit sphere in $\mathbb{R}^p$ and $R_i$ is a nonnegative random variable independent of $u_i$. The covariance matrix can be written as $\bms=p^{-1}E(R_i^2){\bf \Xi}$. Condition (C2) is the same as the conditions (C1) and (C2) in Wang et al. (2015). If the eigenvalues of $\bms$ are all bounded, condition (C2) will hold.

\begin{theorem}\label{th1}
Under Conditions (C1)-(C2), we have $T_S/\sigma_S\cd N(0,1)$ where $\sigma_S^2=\frac{H}{2}\tr^2(\O^2)$.
\end{theorem}

Then, we estimate $\widehat{\tr(\O^2)}$ as
\begin{align*}
\widehat{\tr(\O^2)}=\frac{2}{n(n-1)}\sum_{1\le s<t\le n} (\U_s^\top \U_t)^2.
\end{align*}
By Proposition 1 of \citet{z2023}, we have $\widehat{\tr(\O^2)}/\tr(\O^2)\cp 1$ under the null hypothesis as $n,p\to \infty$. So by Theorem \ref{th1}, we reject the null hypothesis if $T_S/\hat{\sigma}_S>z_{\alpha}$ where $\hat{\sigma}_S^2=\frac{H}{2}\widehat{\tr^2(\O^2)}$ and $z_{\alpha}$ is the upper $\alpha$ quantile of standard normal distribution.

Next, we consider the power function of our test procedure. Specially, we consider the following alternative hypothesis:
\begin{align}\label{h1}
H_1: \bmv_t=\A_0r_t\u_t+\A_1r_{t-1}\u_{t-1}
\end{align}
where $\u_t$ is a random vector uniformly distributed on the unit sphere in $\mathbb{R}^p$ and $r_t$ is a nonnegative random variable independent of $\u_t$.
Let $\bms_0=\A_0^\top \A_0$, $\bms_1=\A_1^\top \A_1$ and $\bms_{01}=\A_0^\top \A_1$. We also assume the following conditions for $\A_0$ and $\A_1$:
\begin{itemize}
\item[(C3)] The eigenvalues of $\bms_0$ are all bounded and $\tr(\bms_1)=O(p/n)$, $\tr(\bms_1^2)=O(p/n)$, $\tr(\bms_{0}\bms_{1})=O(p/n)$.
\end{itemize}

\begin{theorem}\label{th2}
Under $H_1$ in (\ref{h1}) with Condition (C3) holds, if $p/n\to \gamma \in (0,\infty)$, we have, for $H=1$,
\begin{align*}
\frac{T_S-\frac{1}{2}c_1^2\omega^4np^{-2}\tr(\bms_0\bms_1)}{\sqrt{1/2}p^{-2}\omega^4\tr(\bms_0^2)}\cd N(0,1).
\end{align*}
where $c_1=E(r_t)E(r_t^{-1})$ and $\omega=p^{-1/2}\tr^{1/2}(\bms_0)$.
\end{theorem}

Note that under condition (C3), we have $\tr(\O^2)=\frac{\tr(\bms_0^2)}{\tr^2(\bms_0)}(1+o(1))$ under the null hypothesis. So, by Theorem \ref{th2}, the power function of $T_S$ is
\begin{align*}
\beta_S=\lim_{n,p\to \infty}\Phi\left(-z_{\alpha}+\frac{c_1^2n\tr(\bms_0\bms_1)}{\sqrt{2}\tr(\bms_0^2)}\right).
\end{align*}
In addition, according to Theorem 5 in \citet{flm2022}, the power function of the sum-type test proposed by \citet{flm2022} is
\begin{align*}
\beta_{FLM}=\lim_{n,p\to \infty}\Phi\left(-z_{\alpha}+\frac{nE^2(r_t)\tr(\bms_0\bms_1)}{\sqrt{2}E(r_t^2)\tr(\bms_0^2)}\right)
\end{align*}
under Condition (C3). Thus, the asymptotic relative efficiency of our SS test with respect to FLM test is
\begin{align*}
ARE(SS,FLM)=\lim_{p\to \infty}E^2(r_t^{-1})E(r_t^2)\ge\lim_{p\to \infty} \{E(r_t)E(r_t^{-1})\}^2\ge 1
\end{align*}
by Cauchy inequality. Next, we consider three special distributions for $\bmv_t$:
\begin{itemize}
\item[(1)] $\bmv_t\sim N(\bm 0, \I_p)$. So, $r_tE(r_t^{-1}) \cp 1$ and then $ARE(SS,FLM)=1$.
\item[(2)] $\bmv_t\sim t_p(\bm 0, \I_p, v)$. So
\begin{align*}
E(r_t^{-1})=\frac{\Gamma\{(v+1) / 2\}}{v^{1 / 2} \Gamma(v / 2)} \frac{\Gamma\{(p-1) / 2\}}{\Gamma(p / 2)}, E(r_t^2)=\frac{pv}{v-2}
\end{align*}
and then
\begin{align*}
ARE(SS,FLM)=\frac{2}{v-2}\left(\frac{\Gamma((v+1) / 2)}{\Gamma(v / 2)}\right)^2>1.
\end{align*}
For $v=3$, this value is about 2.54 ; for $v=4$, it is about 1.76 ; for $v=\infty$ (multivariate normal distribution), it converges to one.
\item[(3)] $\bmv_t\sim (1-v)N(\bm 0, \I_p)+vN(\bm 0, \sigma^2\I_p)$. So
\begin{align*}
E(r_t^{-1})=\frac{\{v+(1-v) / \sigma\}\left\{v+(1-v) \sigma^2\right\}^{1 / 2}}{2^{1 / 2}} \frac{\Gamma\{(p-1) / 2\}}{\Gamma(p / 2)}, E(r_t^2)=p(1-v+v\sigma^2)
\end{align*}
and then
\begin{align*}
ARE(SS,FLM)=\frac{1+v(1-v)\left(\sigma-\sigma^{-1}\right)^2}{1+v(1-v)\left(1-\sigma^{-1}\right)^2}>1.
\end{align*}
\end{itemize}

\section{Simulations}

We compare our method with the max-type test proposed by Chang, Yao and Zhou (2017) (abbreviated as MAX), the sum-type test (abbreviated as FLM) and the Fisher's combined probability test (abbreviated as FC) proposed by \citet{flm2022}. First, we consider the null hypothesis.  To verify the robustness of the proposed testing method, we consider the following three scenarios for $\bmv_{t}$:
\begin{itemize}
\item[(I)] Multivariate normal distribution. $\bmv_{t}\mathop{\sim}\limits^{i.i.d} N(\boldsymbol  0,\bms)$;
\item[(II)] Multivariate t-distribution. $\bmv_{t} \mathop{\sim}\limits^{i.i.d} t(\boldsymbol  0,\bms,3)$;
\item[(III)] Multivariate mixture normal distribution. $\bmv_{t}$'s are independently generated from  $\gamma f_p(\boldsymbol  0,\bms)+(1-\gamma)f_p(\boldsymbol  0,9\bms)$, denoted by
$\mbox{MN}_{p,\gamma,9}(\boldsymbol  0, \bms)$, where $f_p(\cdot;\cdot)$ is
the density function of $p$-variate multivariate normal
distribution. $\gamma$ is chosen to be 0.8.
\end{itemize}
where $\bms=(\sigma_{ij})_{1\le i,j\le p}$, $\sigma_{ii}=1$, $i=1,\cdots,p$, $\sigma_{ij}=\frac{1}{2}(i-j)^{-2}$ with $i\not=j$. Table \ref{tab:t1} reports the empirical sizes of SS, MAX, FLM, FC tests with $n=100,200$ and $p=40,80,100$. From Table \ref{tab:t1}, we observe that both FLM and SS tests can control the empirical sizes in most cases. However, the empirical sizes of MAX test are a little conservative under the multivariate normal distribution, while a litter larger than the nominal level under the multivariate t-distribution. And the empirical sizes of FC tests also has the same performance as MAX test.

\begin{table}[!th]
\begin{center}
\footnotesize
\caption{\label{tab:t1} Size performance of different tests.}
                     \vspace{0.5cm}
                    \renewcommand\tabcolsep{3.0pt}
                     \renewcommand{\arraystretch}{1.2}
                     {
\begin{tabular}{cc|cccc|cccc|cccc}
\hline \hline
&&\multicolumn{4}{c}{$H=1$}&\multicolumn{4}{c}{$H=2$}&\multicolumn{4}{c}{$H=3$}\\ \hline
$n$&$p$&  MAX &SS &FLM &FC & MAX&SS &FLM &FC & MAX&SS &FLM &FC \\\hline
&&\multicolumn{12}{c}{Multivariate Normal Distribution}\\\hline
100&40&0.011&0.053&0.055&0.042&0.005&0.051&0.052&0.028&0.008&0.042&0.044&0.023\\
100&80&0.008&0.035&0.041&0.022&0.008&0.038&0.039&0.019&0.007&0.059&0.058&0.03\\
100&120&0.003&0.049&0.049&0.029&0.008&0.057&0.059&0.029&0.006&0.053&0.052&0.017\\
200&40&0.015&0.069&0.061&0.046&0.017&0.062&0.058&0.048&0.017&0.037&0.044&0.037\\
200&80&0.014&0.054&0.051&0.038&0.018&0.053&0.059&0.031&0.012&0.05&0.05&0.038\\
200&120&0.01&0.05&0.051&0.035&0.011&0.054&0.05&0.031&0.012&0.044&0.039&0.016\\ \hline
&&\multicolumn{12}{c}{Multivariate t-distribution}\\\hline
100&40&0.068&0.053&0.049&0.079&0.082&0.059&0.052&0.079&0.097&0.065&0.044&0.104\\
100&80&0.067&0.049&0.047&0.074&0.11&0.047&0.054&0.118&0.115&0.039&0.044&0.116\\
100&120&0.088&0.04&0.045&0.091&0.11&0.045&0.036&0.113&0.132&0.045&0.036&0.129\\
200&40&0.075&0.06&0.046&0.095&0.097&0.058&0.056&0.115&0.126&0.061&0.046&0.133\\
200&80&0.1&0.049&0.04&0.107&0.14&0.048&0.029&0.131&0.177&0.061&0.037&0.173\\
200&120&0.113&0.05&0.053&0.11&0.162&0.047&0.048&0.166&0.22&0.06&0.05&0.207\\ \hline
&&\multicolumn{12}{c}{Mixture of Multivariate Normal Distribution}\\\hline
100&40&0.035&0.057&0.058&0.072&0.047&0.068&0.05&0.056&0.029&0.043&0.052&0.051\\
100&80&0.047&0.046&0.048&0.075&0.039&0.041&0.038&0.051&0.043&0.036&0.037&0.045\\
100&120&0.04&0.046&0.045&0.058&0.046&0.048&0.041&0.057&0.033&0.05&0.053&0.049\\
200&40&0.031&0.054&0.055&0.068&0.031&0.052&0.047&0.052&0.025&0.059&0.045&0.049\\
200&80&0.044&0.053&0.061&0.064&0.045&0.052&0.063&0.073&0.038&0.051&0.036&0.066\\
200&120&0.035&0.054&0.049&0.058&0.049&0.049&0.048&0.067&0.045&0.034&0.039&0.064\\
\hline
\hline
\end{tabular}}
\end{center}
\end{table}

Next, we compare the empirical power performance of the above four tests. We consider three models for $\bmv_t$:
\begin{itemize}
\item[(i)] VAR(1) model: $\bmv_t=\A\bmv_{t-1}+\z_t$;
\item[(ii)] VMA(1) model: $\bmv_t=\z_t+\A\z_{t-1}$;
\item[(iii)] VARMA(1) model: $\bmv_t=0.5\A\bmv_{t-1}+\z_t+0.5\A\z_{t-1}$.
\end{itemize}
Here ``VAR(1)", ``VMA(1)" and ``VARMA(1)" are the abbreviations of
1-order vector autoregressive process, vector moving average process and
vector autoregressive moving average process, respectively. Here $\z_t$ are generated from Scenario (I)-(III) with $\bms=\I_p$. Let $\A=(a_{ij})_{1\le i,j\le p}$.  We consider the alternative hypothesis with $a_{ij}\not=0,$ for $1\le i,j\le m$ and $a_{ij}=0$ otherwise. Here $m$ control the signal strength and sparsity of $\A$. We consider two cases for $m$ : (1) dense case: $m=[0.8p]$ and $a_{ij}\sim U(-\frac{1}{4\sqrt{m}},\frac{1}{4\sqrt{m}})$ for $1\le i,j\le m$; (2) sparse case: $m=[0.05p]$ and $a_{ij}\sim U(-\frac{3}{4\sqrt{m}},\frac{3}{4\sqrt{m}})$ for $1\le i,j\le m$.
Table \ref{tab:t2} reports the empirical power of the above four tests with $n=200,p=80$. Under the multivariate normal distribution, the performance of SS test is similar to FLM test, which is consistent to the theoretical result. Under the dense case, the power of sum-type tests--SS and FLM are more powerful than MAX and FC tests. However, MAX and FC tests outperform SS and FLM tests under the sparse case. For the heavy-tailed distributions, our SS test has better performance than FLM test, which shows the advantage of the spatial-sign based method. In addition, we found that the power of the SS and FLM tests with $H=1$ are larger than those tests with $H=2,3$. It is not strange because we consider the alternative hypothesis with 1-order. How to choose the best $H$ for the general case deserves some further studies.

\begin{table}[!th]
\begin{center}
\footnotesize
\caption{\label{tab:t2} Power performance of different tests with $n=200,p=80$.}
                     \vspace{0.5cm}
                    \renewcommand\tabcolsep{3.0pt}
                     \renewcommand{\arraystretch}{1.2}
                     {
\begin{tabular}{cc|cccc|cccc|cccc}
\hline \hline
&&\multicolumn{4}{c}{$H=1$}&\multicolumn{4}{c}{$H=2$}&\multicolumn{4}{c}{$H=3$}\\ \hline
&Model&  MAX &SS &FLM &FC & MAX&SS &FLM &FC & MAX&SS &FLM &FC \\\hline
&&\multicolumn{12}{c}{Multivariate Normal Distribution}\\\hline
Dense&(i)&0.021&0.729&0.725&0.585&0.024&0.497&0.51&0.35&0.02&0.38&0.39&0.262\\
&(ii)&0.016&0.698&0.71&0.559&0.016&0.46&0.465&0.318&0.015&0.358&0.374&0.222\\
&(iii)&0.01&0.732&0.737&0.595&0.016&0.49&0.484&0.348&0.011&0.368&0.379&0.242\\ \hline
Sparse&(i)&0.733&0.613&0.625&0.852&0.701&0.528&0.534&0.802&0.658&0.417&0.423&0.742\\
&(ii)&0.464&0.371&0.374&0.603&0.383&0.234&0.24&0.482&0.326&0.185&0.191&0.411\\
&(iii)&0.603&0.476&0.488&0.738&0.524&0.303&0.32&0.622&0.467&0.254&0.261&0.563\\\hline
&&\multicolumn{12}{c}{Multivariate t-distribution}\\\hline
Dense&(i)&0.193&0.952&0.718&0.673&0.188&0.809&0.482&0.487&0.204&0.637&0.39&0.44\\
&(ii)&0.195&0.939&0.709&0.669&0.192&0.745&0.445&0.456&0.213&0.588&0.366&0.427\\
&(iii)&0.174&0.939&0.736&0.674&0.197&0.758&0.497&0.477&0.238&0.629&0.379&0.449\\\hline
Sparse&(i)&0.766&0.879&0.597&0.862&0.768&0.813&0.506&0.834&0.745&0.708&0.436&0.813\\
&(ii)&0.595&0.559&0.361&0.697&0.546&0.353&0.235&0.609&0.536&0.266&0.194&0.594\\
&(iii)&0.688&0.739&0.456&0.791&0.683&0.567&0.309&0.742&0.636&0.43&0.243&0.684\\\hline
&&\multicolumn{12}{c}{Mixture of Multivariate Normal Distribution}\\\hline
Dense&(i)&0.059&0.931&0.748&0.641&0.048&0.771&0.533&0.429&0.052&0.596&0.379&0.311\\
&(ii)&0.042&0.917&0.723&0.612&0.043&0.703&0.49&0.388&0.058&0.525&0.386&0.305\\
&(iii)&0.042&0.938&0.743&0.609&0.056&0.723&0.47&0.368&0.055&0.571&0.382&0.301\\\hline
Sparse&(i)&0.75&0.846&0.597&0.86&0.736&0.771&0.513&0.807&0.697&0.667&0.444&0.778\\
&(ii)&0.486&0.515&0.363&0.635&0.438&0.341&0.245&0.531&0.399&0.267&0.195&0.472\\
&(iii)&0.635&0.711&0.488&0.778&0.585&0.544&0.323&0.669&0.523&0.413&0.255&0.605\\
\hline
\hline
\end{tabular}}
\end{center}
\end{table}

\section{Appendix}
\subsection{Proof of Theorem \ref{th1}}
Define $$V_{nj}=\sum_{l=1}^H\frac{1}{n-l}\sum_{i=l+1}^{j-1} \U_{i-l}^\top\U_{j-l}\U_i^T \U_j,$$ for $j\in \{3,\cdots,n\}$ and $W_{nk}=\sum_{i=3}^k V_{ni}$, $k\in \{3,\cdots,n\}$. Let $\mathcal{F}_{i}\doteq\sigma\{\U_1,\cdots,\U_i\}$ be the $\sigma$-field generated by $\{\U_j\}_{j\le i}$. It is easy to show that $\mathbb{E}(V_{ni}|\mathcal{F}_{i-1})=0$ and it follows that $\{W_{nk},\mathcal{F}_k: 3\le k\le n\}$ is a zero mean martingale. Let $v_{ni}=\mathbb{E}(V_{ni}^2|\mathcal{F}_{i-1})$, $3\le i \le n$ and $V_n=\sum_{i=3}^n v_{ni}$. The central limit theorem \citep{hh1980} will hold if we can show
\begin{align}\label{clt1}
\frac{V_n}{\var(W_{nn})}\cp 1,
\end{align}
and for any $\epsilon>0$,
\begin{align}\label{clt2}
\sum_{i=3}^n \sigma_S^{-2}
{E}\left[V_{ni}^2\mathbb{I}\left\{|V_{ni}|>\epsilon\sigma_S\right\}|\mathcal{F}_{i-1}\right]
\cp 0.
\end{align}
It can be shown that
\begin{align*}
v_{ni}=&\sum_{h,g=1}^H\frac{1}{(n-h)(n-g)}\sum_{s=h+1}^{i-1}\sum_{t=g+1}^{i-1}E(\U_{s-h}^\top\U_{i-h}\U_s^T \U_i\U_{t-g}^\top\U_{i-g}\U_t^T \U_i\mid \mathcal{F}_{i-1})\\
=&\sum_{h=1}^H\frac{1}{(n-h)^2}\sum_{s=h+1}^{i-1}(\U_{s-h}^\top \U_{i-h})^2\U_s^\top \O\U_s\\
&+\sum_{h=1}^H\frac{2}{(n-h)^2}\sum_{h+1\le s<t\le i-1}\U_{s-h}^\top \U_{i-h}\U_{t-h}^\top \U_{i-h}\U_s^\top \O\U_t
\end{align*}
So
\begin{align*}
\frac{V_n}{\var(W_{nn})}=&\sigma_S^{-2}\sum_{i=3}^n\sum_{h=1}^H\frac{1}{(n-h)^2}\sum_{s=h+1}^{i-1}(\U_{s-h}^\top \U_{i-h})^2\U_s^\top \O\U_s\\
&+\sigma_S^{-2}\sum_{i=3}^n\sum_{h=1}^H\frac{2}{(n-h)^2}\sum_{h+1\le s<t\le i-1}\U_{s-h}^\top \U_{i-h}\U_{t-h}^\top \U_{i-h}\U_s^\top \O\U_t\\
\doteq & C_{n1}+C_{n2}
\end{align*}
Simple algebras lead to
\begin{align*}
E(C_{n1})=\frac{2}{H} \sum_{h=1}^H \frac{1}{(n-h)^2} \sum_{i=h+2}^n(i-h-1)=\frac{1}{H} \sum_{h=1}^H \frac{n-h-1}{n-h} \rightarrow 1,
\end{align*}
as $n\to \infty$. And
\begin{align*}
\var(C_{n1})
& \leq H\sigma_S^{-4} \sum_{h=1}^H \frac{1}{(n-h)^4} \operatorname{Var}\left[\sum_{i=3}^n \sum_{s=h+1}^{i-1} (\U_{s-h}^\top \U_{i-h})^2\U_s^\top \O\U_s\right] \\
& \leq \frac{H}{(n-H)^4\sigma_S^4} \sum_{h=1}^H \operatorname{Var}\left[\sum_{i=3}^n \sum_{s=h+1}^{i-1} (\U_{s-h}^\top \U_{i-h})^2\U_s^\top \O\U_s\right]\\
&=O\left(n^{-1}\frac{E^2((\U_t^\top \O\U_t)^2)}{\tr^4(\O^2)}\right)\to 0
\end{align*}
by Lemma 1 in \citet{wpl}. Thus, we have $C_{n1}\cp 1$. Similarly, $E(C_{n2})=0$ and
\begin{align*}
\var(C_{n2})=O\left\{\frac{E^2\left\{\left(\U_1^T \O \U_2\right)^2\right\}+n^{-1} E^2\left\{\left(\U_1^T \O \U_1\right)^2\right\}}{\tr^4(\O^2)}\right\}\to
0,
\end{align*}
So $C_{n2}\cp 0$. Consequently, (\ref{clt1}) holds.

To show (\ref{clt2}), we only need to prove that
\begin{align*}
\sum_{i=3}^nE(V_{ni}^4)=o(\sigma_S^4).
\end{align*}
By Lemma 1 in \citet{wpl}, we can show that
\begin{align*}
\sum_{i=3}^nE(V_{ni}^4)=O(n^{-2}E^2(\U_1^\top\U_2)^4+n^{-1}E^2((\U_1^\top\U_2)^2(\U_1^\top\U_3)^2))=o(\tr^4(\O^2)).
\end{align*}
Here we complete the proof. \hfill$\Box$

\subsection{Proof of Theorem 2}
\begin{align*}
U(\bmv_t)=&U(\A_0r_t\u_t+\A_1r_{t-1}\u_{t-1})=\frac{\A_0r_t\u_t+\A_1r_{t-1}\u_{t-1}}{||\A_0r_t\u_t+\A_1r_{t-1}\u_{t-1}||}\\
=&\frac{\A_0r_t\u_t+\A_1r_{t-1}\u_{t-1}}{r_tp^{-1/2}\tr^{1/2}(\bms_0)}\frac{r_tp^{-1/2}\tr^{1/2}(\bms_0)}{||\A_0r_t\u_t+\A_1r_{t-1}\u_{t-1}||}\\
=&(\A_0p^{-1/2}\tr^{1/2}(\bms_0)\u_t+\A_1r_{t-1}r_t^{-1}p^{-1/2}\tr^{1/2}(\bms_0)\u_{t-1})(1+\gamma_t)^{-1/2}
\end{align*}
where
\begin{align*}
\gamma_t&=\frac{||\A_0r_t\u_t+\A_1r_{t-1}\u_{t-1}||^2}{r_t^2p^{-1}\tr(\bms_0)}-1\\
&=\frac{r_t^{2}\u_t^\top\bms_0\u_t+2r_tr_{t-1} \u_t^\top\bms_{01}\u_{t-1}+r_{t-1}^2\u_{t-1}^\top\bms_1\u_{t-1} }{r_t^2p^{-1}\tr(\bms_0)}-1\\
&=\left(\frac{\u_t^\top\bms_0\u_t}{p^{-1}\tr(\bms_0)}-1\right)+\frac{2r_{t-1} \u_t^\top\bms_{01}\u_{t-1}}{r_tp^{-1}\tr(\bms_0)}+\frac{r_{t-1}^2\u_{t-1}^\top\bms_1\u_{t-1} }{r_t^2p^{-1}\tr(\bms_0)}\\
&=G_1+G_2+G_3.
\end{align*}
By lemma 4 in \citet{zpfw2014} and condition (C3), we have $E(G_1^2)=O(\tr(\bms_0^2)/p^2)=O(p^{-1}), E(G_2^2)=O(\tr(\bms_{0}\bms_{1})/p^2)=O(p^{-1}n^{-1})$ and $E(G_3^2)=O(p^{-2}(\tr^2(\bms_1)+\tr(\bms_1^2)))=O(n^{-1}p^{-1}+n^{-2})$. So $\gamma_t=O_p(p^{-1/2})$. Thus, by taking the same procedure as the proof of Theorem 1 in \citet{z2022}, we have
\begin{align*}
T_S=&\frac{1}{n-1}\sum_{2\le s<t\le n} \U_{s-1}^\top \U_{t-1}\U_s^\top \U_t\\
=&\frac{1}{n-1}\sum_{2\le s<t\le n}(\A_0p^{-1/2}\tr^{1/2}(\bms_0)\u_{s-1}+\A_1r_{s-2}r_{s-1}^{-1}p^{-1/2}\tr^{1/2}(\bms_0)\u_{s-2})^\top\\
&\times(\A_0p^{-1/2}\tr^{1/2}(\bms_0)\u_{t-1}+\A_1r_{t-2}r_{t-1}^{-1}p^{-1/2}\tr^{1/2}(\bms_0)\u_{t-2})\\
&\times(\A_0p^{-1/2}\tr^{1/2}(\bms_0)\u_s+\A_1r_{s-1}r_s^{-1}p^{-1/2}\tr^{1/2}(\bms_0)\u_{s-1})^\top\\
&\times(\A_0p^{-1/2}\tr^{1/2}(\bms_0)\u_t+\A_1r_{t-1}r_t^{-1}p^{-1/2}\tr^{1/2}(\bms_0)\u_{t-1})+o_p(p^{-1})
\end{align*}
Let $\omega=p^{-1/2}\tr^{1/2}(\bms_0), \delta_t=r_{t-1}r_t^{-1}$. We can decompose $T_S$ as
\begin{align*}
T_S=&\frac{\omega^4}{n-1} \sum_{2\le s < t\le n}\boldsymbol{u}_s^{\top} \mathbf{A}_0^{\top} \mathbf{A}_0 \boldsymbol{u}_t \boldsymbol{u}_{t-1}^{\top} \mathbf{A}_0^{\top} \mathbf{A}_0 \boldsymbol{u}_{s-1}+\frac{\omega^4}{n-1} \sum_{2\le s < t\le n}\delta_{t-1}\delta_{s-1}\boldsymbol{u}_{s-1}^{\top} \mathbf{A}_1^{\top} \mathbf{A}_1 \boldsymbol{u}_{t-1} \boldsymbol{u}_{t-1}^{\top} \mathbf{A}_0^{\top} \mathbf{A}_0 \boldsymbol{u}_{s-1}\\
&+D_1+D_2+D_3+o_p(p^{-1})
\end{align*}
where
\begin{align*}
D_1=
& \frac{\omega^4}{n-1} \sum_{2\le s < t\le n}\left(\delta_{s-1}\delta_{s-2}\delta_{t-1}\delta_{t-2}\boldsymbol{u}_{s-1}^{\top} \mathbf{A}_1^{\top} \mathbf{A}_1 \boldsymbol{u}_{t-1} \boldsymbol{u}_{t-2}^{\top} \mathbf{A}_1^{\top} \mathbf{A}_1 \boldsymbol{u}_{s-2}+\delta_{s-2}\delta_{t-2}\boldsymbol{u}_s^{\top} \mathbf{A}_0^{\top} \mathbf{A}_0 \boldsymbol{u}_t \boldsymbol{u}_{t-2}^{\top} \mathbf{A}_1^{\top} \mathbf{A}_1 \boldsymbol{u}_{s-2}\right)\\
D_2=
& \frac{\omega^4}{n-1} \sum_{2\le s < t\le n}\left(\delta_{t-1}\boldsymbol{u}_s^{\top} \mathbf{A}_0^{\top} \mathbf{A}_1 \boldsymbol{u}_{t-1} \boldsymbol{u}_{t-1}^{\top} \mathbf{A}_0^{\top} \mathbf{A}_0 \boldsymbol{u}_{s-1}+\delta_{s-1}\boldsymbol{u}_{s-1}^{\top} \mathbf{A}_1^{\top} \mathbf{A}_0 \boldsymbol{u}_t \boldsymbol{u}_{t-1}^{\top} \mathbf{A}_0^{\top} \mathbf{A}_0 \boldsymbol{u}_{s-1}\right. \\
& +\delta_{s-1}\delta_{t-1}\boldsymbol{u}_{s-1}^{\top} \mathbf{A}_1^{\top} \mathbf{A}_1 \boldsymbol{u}_{t-1} \boldsymbol{u}_{t-1}^{\top} \mathbf{A}_0^{\top} \mathbf{A}_1 \boldsymbol{u}_{s-2}+\delta_{s-1}\delta_{t-1}\delta_{t-2}\boldsymbol{u}_{s-1}^{\top} \mathbf{A}_1^{\top} \mathbf{A}_1 \boldsymbol{u}_{t-1} \boldsymbol{u}_{t-2}^{\top} \mathbf{A}_1^{\top} \mathbf{A}_0 \boldsymbol{u}_{s-1} \\
& +\delta_{t-2}\boldsymbol{u}_s^{\top} \mathbf{A}_0^{\top} \mathbf{A}_0 \boldsymbol{u}_t \boldsymbol{u}_{t-2}^{\top} \mathbf{A}_1^{\top} \mathbf{A}_0 \boldsymbol{u}_{s-1}+\delta_{s-2}\boldsymbol{u}_s^{\top} \mathbf{A}_0^{\top} \mathbf{A}_0 \boldsymbol{u}_t \boldsymbol{u}_{t-1}^{\top} \mathbf{A}_0^{\top} \mathbf{A}_1 \boldsymbol{u}_{s-2}\\
&\left.+\delta_{t-1}\delta_{t-2}\delta_{s-2}\boldsymbol{u}_s^{\top} \mathbf{A}_0^{\top} \mathbf{A}_1 \boldsymbol{u}_{t-1} \boldsymbol{u}_{t-2}^{\top} \mathbf{A}_1^{\top} \mathbf{A}_1 \boldsymbol{u}_{s-2}+\delta_{s-1}\delta_{t-2}\delta_{s-2}\boldsymbol{u}_{s-1}^{\top} \mathbf{A}_1^{\top} \mathbf{A}_0 \boldsymbol{u}_t \boldsymbol{u}_{t-2}^{\top} \mathbf{A}_1^{\top} \mathbf{A}_1 \boldsymbol{u}_{s-2}\right)\\
D_3=
& \frac{\omega^4}{n-1} \sum_{2\le s < t\le n}\left(\delta_{t-1}\delta_{s-2}\boldsymbol{u}_s^{\top} \mathbf{A}_0^{\top} \mathbf{A}_1 \boldsymbol{u}_{t-1} \boldsymbol{u}_{t-1}^{\top} \mathbf{A}_0^{\top} \mathbf{A}_1 \boldsymbol{u}_{s-2}+\delta_{s-1}\delta_{t-2}\boldsymbol{u}_{s-1}^{\top} \mathbf{A}_1^{\top} \mathbf{A}_0 \boldsymbol{u}_t \boldsymbol{u}_{t-2}^{\top} \mathbf{A}_1^{\top} \mathbf{A}_0 \boldsymbol{u}_{s-1}\right. \\
& \left.+\delta_{t-1}\delta_{t-2}\boldsymbol{u}_s^{\top} \mathbf{A}_0^{\top} \mathbf{A}_1 \boldsymbol{u}_{t-1} \boldsymbol{u}_{t-2}^{\top} \mathbf{A}_1^{\top} \mathbf{A}_0 \boldsymbol{u}_{s-1}+\delta_{s-1}\delta_{s-2}\boldsymbol{u}_{s-1}^{\top} \mathbf{A}_1^{\top} \mathbf{A}_0 \boldsymbol{u}_t \boldsymbol{u}_{t-1}^{\top} \mathbf{A}_0^{\top} \mathbf{A}_1 \boldsymbol{u}_{s-2}\right)
\end{align*}
After some tedious algebra, we have
\begin{align*}
E(D_1^2)=o(p^{-2}), E(D_2^2)=o(p^{-2}), E(D_3^2)=o(p^{-2})
\end{align*}
by Condition (C3). Taking the same procedure as the proof of Theorem \ref{th1}, we can show that
\begin{align*}
\frac{1}{(n-1)\sqrt{p^{-4}\tr^2(\bms_0^2)/2}} \sum_{2\le s < t\le n}\boldsymbol{u}_s^{\top} \mathbf{A}_0^{\top} \mathbf{A}_0 \boldsymbol{u}_t \boldsymbol{u}_{t-1}^{\top} \mathbf{A}_0^{\top} \mathbf{A}_0 \boldsymbol{u}_{s-1}\cd N(0,1).
\end{align*}
And
\begin{align*}
&E\left(\frac{\omega^4}{n-1} \sum_{2\le s < t\le n}\delta_{t-1}\delta_{s-1}\boldsymbol{u}_{s-1}^{\top} \mathbf{A}_1^{\top} \mathbf{A}_1 \boldsymbol{u}_{t-1} \boldsymbol{u}_{t-1}^{\top} \mathbf{A}_0^{\top} \mathbf{A}_0 \boldsymbol{u}_{s-1}\right)=\frac{1}{2}c_1^2\omega^4np^{-2}\tr(\bms_0\bms_1),\\
&\var\left(\frac{\omega^4}{n-1} \sum_{2\le s < t\le n}\delta_{t-1}\delta_{s-1}\boldsymbol{u}_{s-1}^{\top} \mathbf{A}_1^{\top} \mathbf{A}_1 \boldsymbol{u}_{t-1} \boldsymbol{u}_{t-1}^{\top} \mathbf{A}_0^{\top} \mathbf{A}_0 \boldsymbol{u}_{s-1}\right)=c_1^2\omega^4np^{-2}\tr(\bms_0\bms_1)=o(p^{-2}).
\end{align*}
Thus, we have
\begin{align*}
\frac{T_S-\frac{1}{2}c_1^2\omega^4np^{-2}\tr(\bms_0\bms_1)}{\sqrt{1/2}p^{-2}\omega^4\tr(\bms_0^2)}\cd N(0,1).
\end{align*}
\hfill $\Box$


\end{document}